\documentclass[pamm,a4paper,fleqn]{w-art}

\usepackage{times,cite,w-thm}
\usepackage[T1]{fontenc}
\usepackage[utf8]{inputenc}
\usepackage{graphicx}
\usepackage{nicefrac}

\usepackage{amsmath}
\usepackage{amssymb}

\newcommand{\surf}{\mathcal{S}}

\newcommand{\R}{\mathbb{R}}

\newcommand{\vnor}{\mathcal{V}}
\newcommand{\normal}{\boldsymbol{\nu}}

\newcommand{\vel}{\boldsymbol{u}}

\newcommand{\velTest}{\boldsymbol{v}}
\newcommand{\relvel}{\boldsymbol{w}}

\newcommand{\param}{\boldsymbol{X}}

\newcommand{\paramTest}{\boldsymbol{Z}}

\newcommand{\paramUpdate}{\boldsymbol{Y}}

\newcommand{\coord}{\boldsymbol{x}}

\newcommand{\shop}{\boldsymbol{B}}

\newcommand{\meanc}{\mathcal{H}}

\newcommand{\proj}{\boldsymbol{P}}
\newcommand{\Id}{\boldsymbol{I}}

\newcommand{\DivP}{\operatorname{div}_{\!\proj}\!}
\newcommand{\DivC}{\operatorname{div}_{\!C}\!}

\newcommand{\inputTikzPic}[1]{\ifthenelse{\boolean{plotTikzPics}}{\input{#1}}{\fbox{\centering\begin{minipage}[t][5cm][t]{0.9\textwidth}\centering\textbf{\color{red}enable plotTikzPics}\end{minipage}}}}

\newcommand{\Grad}{\nabla}
\newcommand{\Div}{\operatorname{div}}%
\newcommand{\DivSurf}{\Div_{\!\surf}}%
\newcommand{\GradSurf}{\Grad_{\!\surf}}
\newcommand{\GradP}{\Grad_{\!\proj}}
\newcommand{\GradC}{\Grad_{\!C}}

\newcommand{\ReynoldsNumber}{\mathrm{Re}}

\renewcommand{\Re}{\ReynoldsNumber}

\newcommand{\stressC}{\sigma}
\newcommand{\stress}{\boldsymbol{\stressC}}

\begin{document}

\TitleLanguage[EN]
\title[The short title]{The title}

\title[]{A surface finite element method for the Navier-Stokes equations on evolving surfaces}

\author{\firstname{Veit} \lastname{Krause}\inst{1}}
\author{\firstname{Eric} \lastname{Kunze}\inst{1}}
\author{\firstname{Axel} \lastname{Voigt}\inst{1,2,3}}

\address[\inst{1}]{\CountryCode[DE]Institute of Scientific Computing, TU Dresden, Dresden, 01062, Germany}

\address[\inst{2}]{\CountryCode[DE]Center for Systems Biology Dresden (CSBD), Pfotenhauerstra{\ss}e 108, Dresden, 01307, Germany}

\address[\inst{3}]{\CountryCode[DE]Cluster of Excellence Physics of Life (PoL), Dresden, 01062, Germany}

\AbstractLanguage[EN]
\begin{abstract}
We introduce a surface finite element method for the numerical solution of Navier–Stokes equations on evolving surfaces with a prescribed deformation of the surface in normal direction. The method is based on approaches for the full surface Navier-Stokes equations in the context of fluid-deformable surfaces and adds a penalization of the normal component. Numerical results demonstrate the same optimal order as proposed for surface (Navier-)Stokes equations on stationary surfaces. The approach is applied to high-resolution 3D scans of clothed bodies in motion to provide interactive virtual fluid-like clothing.  
\end{abstract}
\maketitle                   

\section{Introduction}

\begin{figure}
    \centering
    \includegraphics[width=\linewidth]{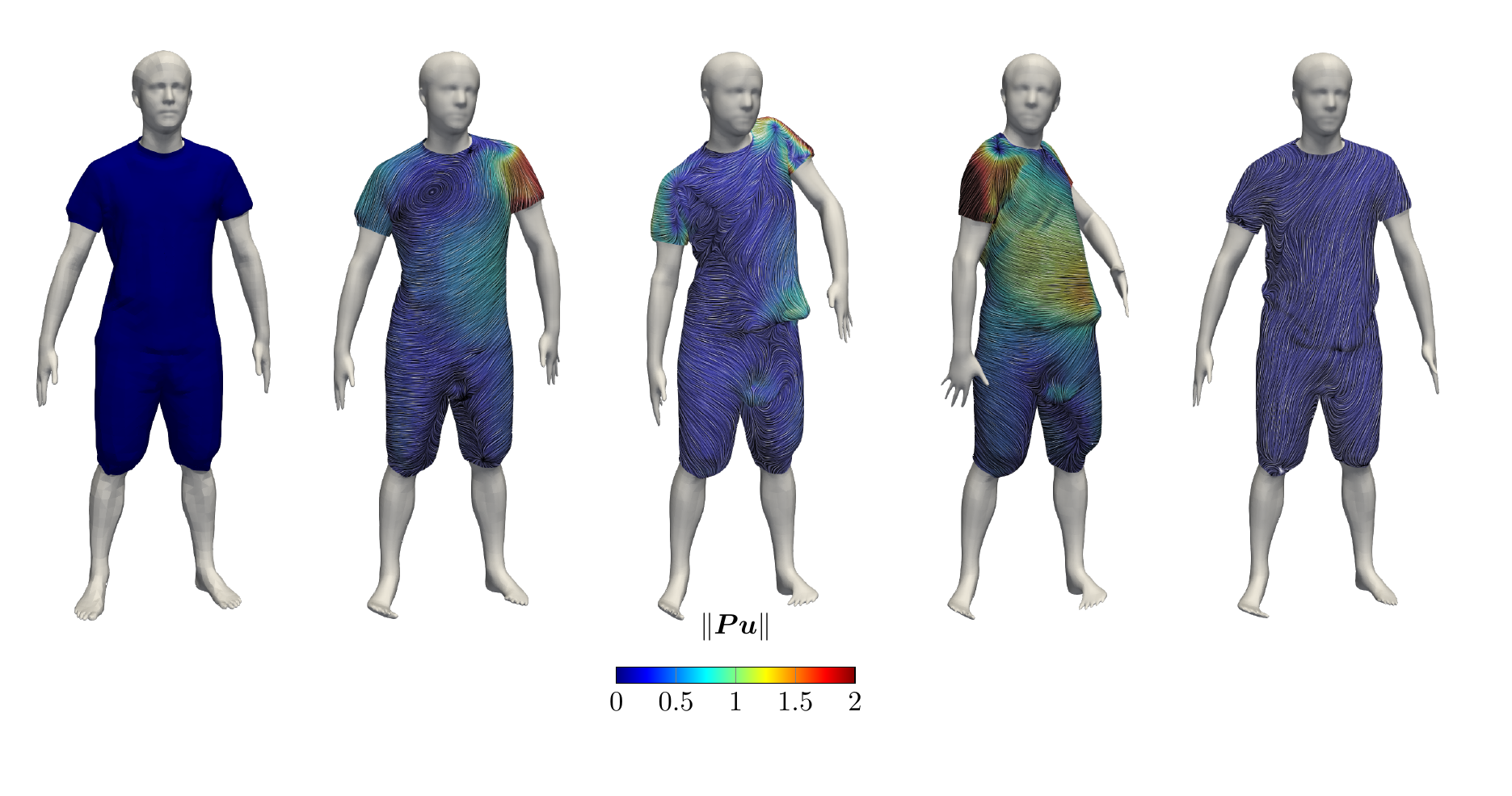}
    \caption{Snapshots of a clothed body in motion with virtual fluid-like clothing. The fluid flow is induced by the motion and visualized by a LIC filter with color coding the magnitude of the tangential velocity.} 
    \label{fig:clothing}
\end{figure}

Fluid deformable surfaces are ubiquitous interfaces in biology, playing an
essential role in processes from the subcellular to the tissue scale. For instance,
lipid membranes are fluidic thin sheets that define the boundaries of cells and
compartmentalize them. From a mechanical point of view, lipid membranes are
soft materials exhibiting a solid–fluid duality: while they store elastic energy when
stretched or bent, like solid shells, under in-plane shear they flow as two-dimensional, viscous fluids. With this solid–fluid duality any shape change contributes to tangential flow and vice versa any tangential flow on a curved surface induces shape deformations. The mathematical description of fluid deformable surfaces has been introduced in \cite{Arroyoetal_PRE_2009,Salbreuxetal_PRE_2017,torres2019modelling,reuther2020numerical} and is numerically solved in \cite{torres2019modelling,reuther2020numerical,Krause22,bachini2023derivation}. We here consider a simplified model with a prescribed deformation of the surface. What remains to solve is the fluid part, the Navier-Stokes equations on evolving surfaces. Two different numerical approaches exist for this problem, a surface finite element method for a vorticity-stream function formulation \cite{Reutheretal_MMS_2015,Reutheretal_MMS_2018} and surface or trace finite element methods in the principal variables \cite{Nitschkeetal_PRF_2019,olshanskii2022tangential,olshanskii2023}. Both approaches are based on a splitting of the velocity in a tangential and a normal component and solely address the resulting tangential Navier-Stokes equations. We here propose a different approach, which is based on the full problem of the velocity and penalization of the normal component to the prescribed velocity. Besides the numerical features of this approach we demonstrate the applicability of the approach for realistic virtual clothing, see figure \ref{fig:clothing}. High-resolution 3D scans of clothed bodies in motion from the CAPE project \cite{Pons-Moll:Siggraph2017,CAPE:CVPR:20} provide the raw data for the prescribed normal velocity. This motion induces tangential flow which is visualized to create virtual interactive clothing. For a detailed explanation of the pre- and postprocessing, we refer to \cite{kunze2022}.  

The paper is structured as follows. We introduce the notation and the surface operators and explain the integration of the constraint into the full surface Navier-Stokes equations in section \ref{sec:SNS}. In section \ref{sec:numeric} we explain an approach for area conservative surface movement with mesh regularization and discretization of the model by surface finite elements (SFEM). Based on this discretization we show numerical convergence studies and simulation results in section \ref{sec:results}. In section \ref{sec:conclusions} we draw conclusions.

\section{Surface Navier-Stokes equations on evolving surfaces}
\label{sec:SNS}

We consider a time-dependent smooth domain $\Omega=\Omega(t)\subset\R^3$ and its boundary $\surf = \surf(t)$, which is given by a parametrization $\param$. Related to $\surf$ we denote the surface normal $\normal$, the shape Operator $\shop= -\GradC \normal$, the mean curvature $\meanc= \operatorname{tr}\shop$ and the surface projection $\proj=\Id-\normal\otimes\normal$. 
Let $\vel$ be a continuously differentiable $\R^3$-vector field and $\stress$ a continuously differentiable $\R^{3\times3}$-tensor field defined on $\surf$. 
As in \cite{jankuhn2018,Nitschke2022GoI,bachini2023derivation} we define the surface tangential gradient $\GradP$ and the componentwise surface gradient $\GradC$ by
$\GradP\vel = \proj\nabla\vel^e\proj$ and
$\GradC\stress = \nabla\stress^e \proj$, where $\vel^e$ and $\stress^e$ are arbitrary smooth extensions of $\vel$ and $\stress$ in the normal direction and $\nabla$ is the gradient of the embedding space $\R^3$. We define the corresponding divergence operators by $\DivP\vel=\operatorname{tr}(\GradP\vel)$ and $\DivC(\stress\proj)= \operatorname{tr}\GradC(\stress\proj)$ where $\operatorname{tr}$ is the trace operator. The last leads to a non-tangential vector field even if $\stress$ is a tangential tensor field. The operator $\DivP$ relates to the covariant divergence by $\DivP\vel = \DivSurf\proj\vel-\vel\normal\meanc$, see \cite{Nitschke2022GoI}. 

On the surface $\surf$ we consider the material velocity $\vel$. Similar to \cite{reuther2020numerical,Krause22,bachini2023derivation} we use a Lagrangian description in the normal direction and Eulerian description in tangential direction where the material relates to the surface by $\partial_t\param\normal =\vel\normal$. We write the material derivative in the form $\dot{\vel}=\partial_t\vel+\nabla_{\relvel}\vel$ where $\relvel=\vel-\partial_t\param$ is the relative material velocity \cite{Nitschke2022GoI,NitschkeVoigt2023Tensorial}. We define the viscous stress tensor $\stress=\frac{1}{2}(\GradP\vel+\GradP\vel^T)$ as in \cite{jankuhn2018,Brandneretal_SIAMJSC_2022,torres2019modelling,Krause22,bachini2023derivation} and consider the full surface Navier-Stokes equations \cite{Krause22,bachini2023derivation,reuther2020numerical,jankuhn2018}
\begin{align}
    \label{eq:navierstokes}
    \begin{aligned}
        \partial_t \vel + \nabla_{\relvel}\vel &= -\GradSurf p-p\meanc\normal + \frac{2}{\Re}\DivC \stress + \mathbf{b}  \\
        \DivP \vel &= 0  
    \end{aligned}    
\end{align}
with pressure $p$, possible forcing term $\mathbf{b}$ and Reynolds number $\Re$. In the following, we assume that the normal velocity of the surface $\vnor$ is given. This lead to the constraint $\vel \cdot \normal = \vnor$. However, $\vnor$ cannot be arbitrary. With the identities described above the local inextensibility constraint $\DivP\vel=0$ implies
\begin{align}
  \label{eq:vnorconst}
  0 = \int_{\surf} \DivP\vel = \int_{\surf} \DivSurf\proj\vel - \vnor\meanc = -\int_{\surf}\vnor\meanc .    
\end{align}
This leads to the condition $\vnor\in \{\vnor\in H^1(\surf) \vert \int_\surf \vnor\meanc=0 \}$. In other words, the given normal velocity $\vnor$ has to fulfill global area conservation.
We add the constraint to the formulation of the surface Navier-Stokes equations \eqref{eq:navierstokes} by a $L^2$-penalisation term leading to 
\begin{align}
    \label{eq:navierstokesIII}
    \begin{aligned}
        \partial_t \vel + \nabla_{\relvel}\vel &= -\GradSurf p-p\meanc\normal + \frac{2}{\Re}\DivC \stress + \mathbf{b}  + \beta(\vel\normal-\vnor)\normal   \\
        \DivP \vel &= 0  
    \end{aligned}    
\end{align}
with $\beta>0$ a penalisation parameter. The approach has similarities with previous approaches to solve the surface Navier-Stokes equations or other vector-valued surface partial differential equations on stationary surfaces  \cite{Nestleretal_JNS_2018,jankuhn2018,Nitschkeetal_PRSA_2018,Reutheretal_PF_2018,nestler2019finite,reuther2020numerical,Hansboetal_IMAJNA_2020,Brandneretal_SIAMJSC_2022,Bachinietal_arXiv_2022}. In these approaches penalization of the normal component is used to ensure tangentiality of a $\R^3$-vector field $\mathbf{u}$ by $\beta(\vel \cdot \normal - 0)\normal$. A general error analysis for this approach is done in \cite{hardering2021tangential}. Here we extend this concept for the penalization to a given scalar field and assume to obtain the same convergence results as in \cite{hardering2021tangential}. Compared to \cite{Reutheretal_MMS_2015,Reutheretal_MMS_2018,Nitschkeetal_PRF_2019,olshanskii2022tangential}, where an equation for the tangential velocity is solved, a solution of eq. \eqref{eq:navierstokesIII} fulfills the normal constraint only up to penalization error. We consider an experimental convergence study of this error in section \ref{sec:results}. 

\section{Numerical approach}
\label{sec:numeric}
\subsection{Surface representation}
In the following, we assume a piecewise smooth approximation $\surf_h$ of the surface $\surf$.
Let $\surf_h^{lin}$ be a piecewise linear reference surface given by shape regular triangulation $\mathcal{T}^{\mathrm{lin}}$ and $\param$ a bijective map $\param \colon \surf_h^{lin} \rightarrow \surf$ such that $\surf=\bigcup_{\hat{T}\in \mathcal{T}^{\mathrm{lin}}} \param(\hat{T})$. The construction of such maps is discussed in \cite{praetorius2020dunecurvedgrid,Brandneretal_SIAMJSC_2022}. 

By the $k$-th order interpolation $\param_k=I^k(\param)$ we get a piecewise polynomial surface $\surf_h = \bigcup_{\hat{T}\in \mathcal{T}^{\mathrm{lin}}} \param_k(\hat{T})$. We use each geometrical quantity, like the normal vector $\normal_h$, the shape operator $\shop_h$ and the inner product $(\cdot, \cdot)_h$, with respect to $\surf_h$. We denote the size of the grid by $h$, i.e. the longest edge of the mesh. We discretize the surface in time at timesteps $t^n$ by $\surf_h^n = \surf_h(t^n)$ where the normal vector $\normal_h$ is chosen with respect $S_h^n$ and we define a discrete function space $V_k(\surf_h^n)=\{ v \in C^0(\surf_h) \vert \quad v\vert_{T} \in \mathbb{P}_k(T) \}$. 

\subsection{Surface evolution with area conservation}
To test our approach we need to specify a surface evolution that fulfills the global area constraint. For a given $\vnor_0(\meanc)\in H^1$ we can construct a surface evolution with global area constraint by
\begin{align}
    \label{eq:areaconst}
    \vnor(\meanc) = \vnor_0(\meanc) - \frac{\int_\surf \vnor_0(\meanc)\meanc}{\int_\surf \meanc}
\end{align}    
as discussed in \cite{barrett2008parametric}. The surface can than be moved by $\partial_t\param=\vnor(\meanc)\normal$ a pure normal movement. Alternatively, the surface evolution can be combined with a numerical approach for mesh regularization as presented in \cite{barrett2008parametric}
\begin{align}
  \label{eq:meshmovement}
  \begin{aligned}
      \partial_t \param  \normal &= \vnor(\meanc)\\
      \meanc\normal &= \Delta_C \param. 
  \end{aligned}    
\end{align}
In addition to the normal movement of the surface, we get a tangential movement which leads to mesh regularization after surface discretization. Moreover, we get a representation of the mean curvature $\meanc$ which is discussed by \cite{Dziuk_NM_1990,Baenschetal_JCP_2005,Hausseretal_JSC_2007,Dziuk_NM_2008} and it can be used for a numerically stable time discretization. The quality of this approach is discussed for mean curvature flow $\vnor = \meanc$ and surface diffusion $\vnor=\Delta_{\surf}\meanc$ in \cite{barrett2008parametric} and used for fluid deformable surfaces in \cite{Krause22,bachini2023derivation} by $\vnor=\vel \cdot \normal$. Here this approach is used for the global area conserving velocity $\vnor$. 
At each timestep $t^n=n\tau$ with $\tau>0$ we solve eqs. \eqref{eq:meshmovement} and \eqref{eq:areaconst}. To do so we insert the surface update $\paramUpdate^{n+1}=\param_k^{n+1}-\param_k^n$ in eq. \eqref{eq:meshmovement} to find $\paramUpdate^{n+1}\in V_k(\surf_h^n)^3$ and $\meanc^{n+1}\in V_k(\surf_h^n)$. The nonlinear system is solved by an iteration over $j>0$ of the linear system 
\begin{align}
  \label{eq:meshmovementII}
  \begin{aligned}
      (\paramUpdate^{j+1}\normal_h,w) &= \left(\tau\vnor_0(\meanc^{j+1}) - \frac{\tau\int_{\surf^j} \paramUpdate^j\normal_h\meanc^j}{\int_{\surf^j} \meanc^j} , w\right)\\
      (\meanc^{j+1}\normal_h , \paramTest)_h + (\GradC\paramUpdate^{j+1} ,  \GradC\paramTest)_h &= (\GradC\param^n ,  \GradC\paramTest)_h 
  \end{aligned}    
\end{align}
with appropriate test functions $w$ and $\mathbf{Z}$ until the error $\int_{\surf^j} \paramUpdate^j\normal_h\meanc^j/ \int_{\surf^j} \meanc^j <\epsilon$, with $\epsilon > 0$ a small tolerance. After each iteration we set $\paramUpdate^{n+1} = \paramUpdate^{j+1}$ and $\meanc^{n+1} = \meanc^{j+1}$ and compute the new surface $\surf^{n+1}$. 

\subsection{Discretization of the surface Navier-Stokes equations}
The numerical approach closely follows \cite{Krause22}. However, as $\surf^{n+1}$ is known we can discretize the surface Navier-Stokes equations \eqref{eq:navierstokesIII} in each time step with respect to the new surface $\surf^{n+1}$ to find the new velocity $\vel^{n+1} \in V_k(\surf_h^{n+1})^3$ and new pressure $p^{n+1} \in V_{k-1}(\surf_h^{n+1})$. We consider this Taylor-Hood element with $k = 3$. In each timestep, we solve the linear equations  
\begin{align}
    \label{eq:navierstokesIV}
    \begin{aligned}
        (\vel^{n+1}+\tau\nabla_{\relvel^n}\vel^{n+1},\velTest)_h =& (\hat{\vel}^n,\velTest)_h
        +(\tau p^{n+1},\DivP\velTest)_h 
        -\frac{2\tau}{\Re}(\stress^{n+1},\GradP\velTest)_h 
        +(\tau\mathbf{b},\velTest)_h  \\
        &+\beta^n\tau(\vel^{n+1}\normal_h-\vnor^{n+1} , \velTest\normal_h)_h   \\
        (\DivP \vel^{n+1},q)_h =& 0  
    \end{aligned}    
\end{align} 
where $\relvel^n = \hat{\vel}^n - \frac{1}{\tau}\paramUpdate^{n+1}$, with $\hat{\vel}^n$ the lift of $\vel^n$ defined on $\surf^n$ to $\surf^{n+1}$. We use the penalization parameter $\beta^n = \beta_0 h^{-2}$ with $\beta_0=100$. The scaling with mesh size $h$ guarantees optimal order of convergence in the corresponding case for stationary surfaces \cite{hardering2021tangential}. Compared to the computations in \cite{Krause22,bachini2023derivation} here the fluid velocity fulfills the inextensibility constraint $\DivP\vel^{n+1}=0$ and the normal constraint $\vel^{n+1} \cdot \normal_h=\vnor^{n+1}$ with respect to the new surface $\surf^{n+1}$ up to numerical errors. 

\subsection{Comments on the implementation}
Implementation is done by the finite element framework AMDiS \cite{vey2007amdis,witkowski2015software} base on DUNE \cite{SanderDune2020,bastian2021dune} using the grid manager AluGrid \cite{dunealugrid:16} and the CurvedGrid library \cite{praetorius2020dunecurvedgrid}. All linear systems are solved by a direct solver using the external library PETSc.  

\section{Results}
\label{sec:results}
\subsection{Flow on a perturbed sphere}
\begin{figure}
    \centering
    \includegraphics[width=\linewidth]{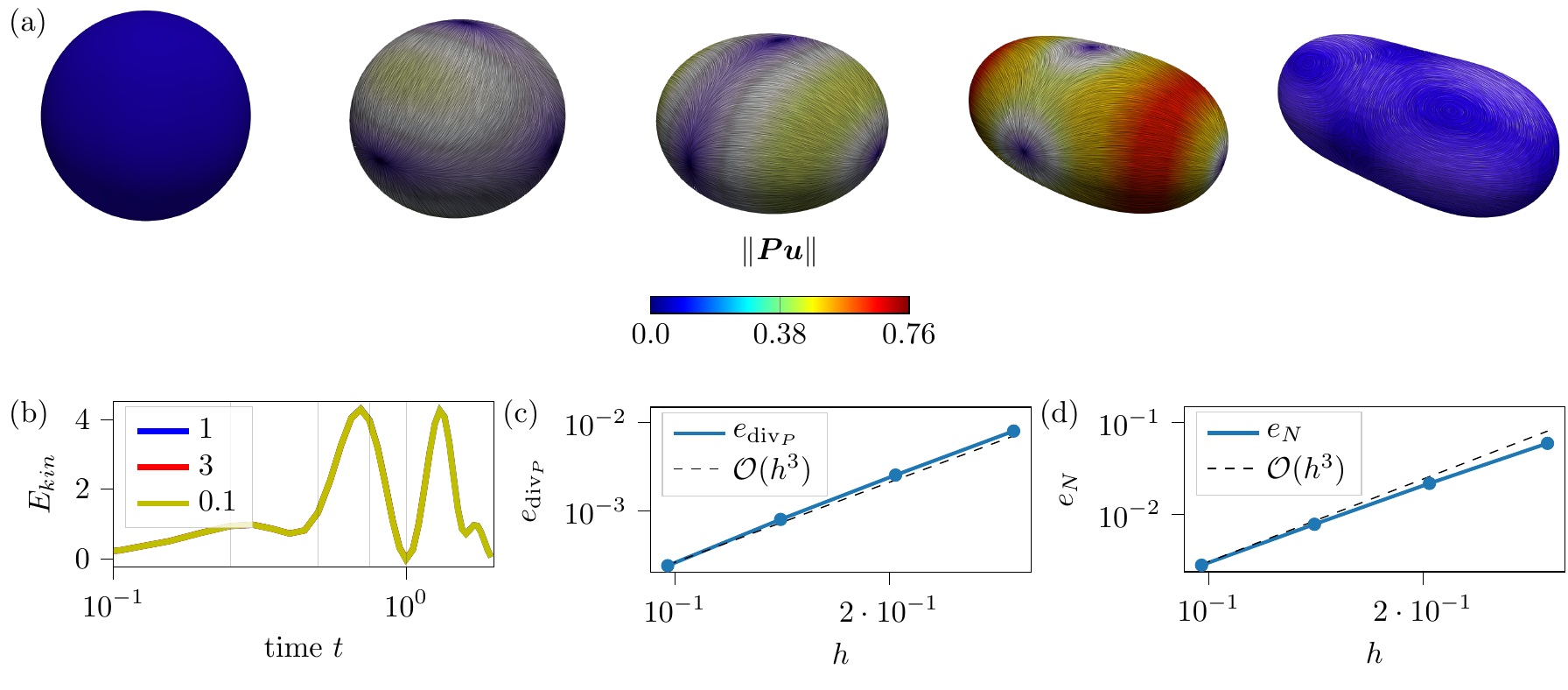}
    \caption{(a) Snapshots of the surface at $t=0,0.25,0.5,0.75,1$. Visualized is the tangential velocity $\proj\vel$ by a LIC Filter and its magnitude by the color coding. (b) Kinetic energy $E_{kin}$ over time for different Reynolds numbers. The time steps shown in (a) are marked. (c,d) The convergence study for the incompressible error and normal error.}
    \label{fig:vis1}
\end{figure}

\begin{figure}
    \centering
    \includegraphics[width=\linewidth]{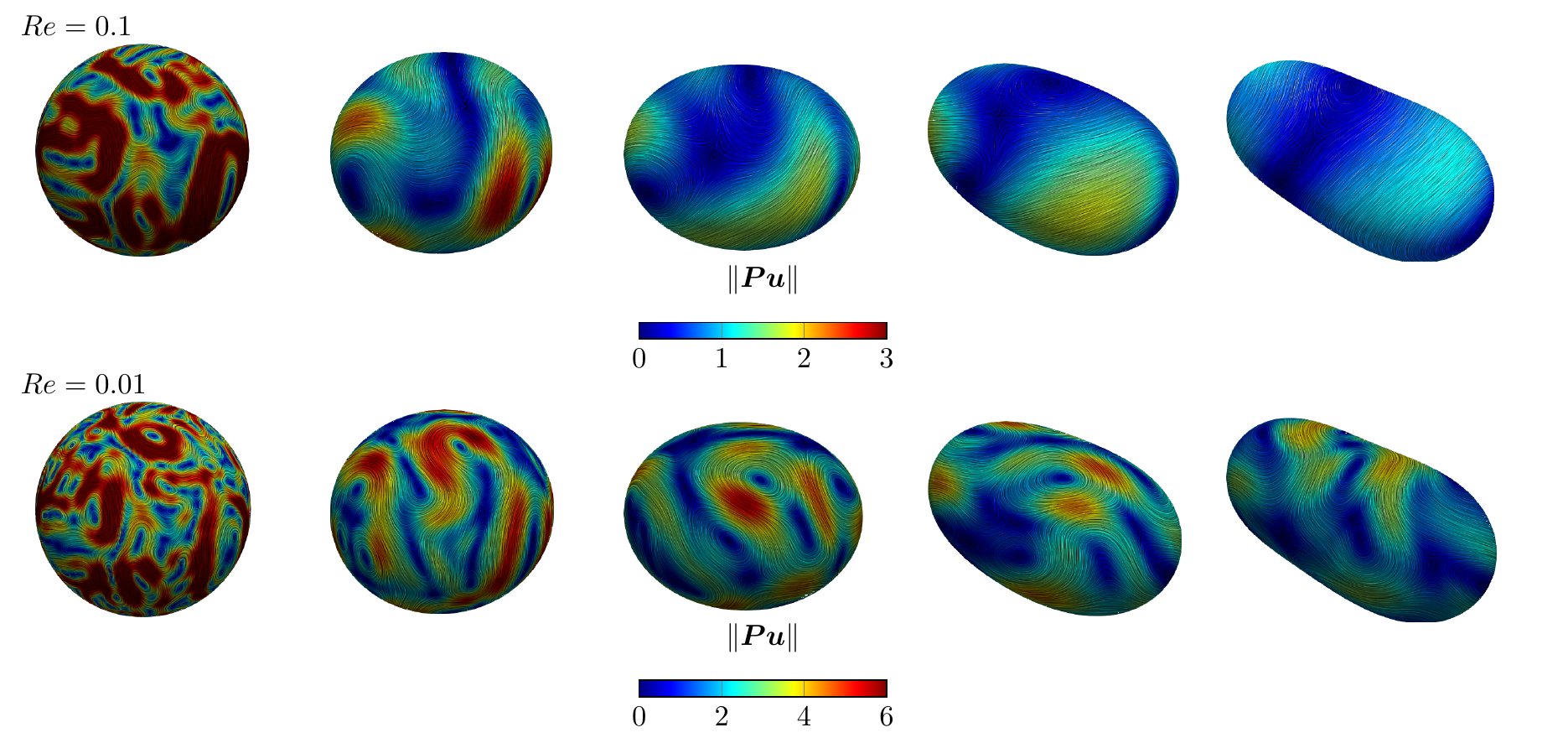}
    \caption{Snapshots of the surface at $t=0,0.25,0.5,0.75,1$ initialized with random initial value for different Reynolds numbers. Visualized is the tangential flow $\proj\vel$ by a LIC Filter and its magnitude by the color coding.}
    \label{fig:vis2}
\end{figure}

We initialize the surface $\surf$ by the parametrization of the unit sphere and define the given normal velocity by $\vnor_0(\meanc) = \alpha\meanc + \sin(\pi t) x_0^2 + \sin(2\pi t)x^2_1$
with $\alpha=1e-3$, which is time-periodic on the interval $[0,2)$. The polynomial order for the surface and the Taylor-Hood element is $k=3$. The evolution of the surface is done by eqs. \eqref{eq:meshmovementII} and the Navier-Stokes equation by eqs. \eqref{eq:navierstokesIV}. 
In figure \ref{fig:vis1} (a) the evolution of the surface is shown for the initial fluid velocity $\vel_0=0$. 
In (b) we show the kinetic energy $E_{kin}=\int_{\surf} \Vert\vel\Vert^2$ for different Reynolds numbers $\Re$. The energy is independent of $\Re$ which indicates that the induced flow is independent too. By following the stream-function formulation in \cite{Reutheretal_MMS_2015,Reutheretal_MMS_2018} the solution of the surface Navier-Stokes equations \eqref{eq:navierstokesIII} can be decomposed orthogonally into a $\DivSurf$-free and $\DivP$-free component where the $\DivP$-free component is defined by a Laplace equation and independent from the initial value on the tangential flow or $\Re$ in case of a given normal velocity.  
In figure \ref{fig:vis2} the same simulation is done for a random initial value for the fluid flow. Compared to the results in figure \ref{fig:vis1} there is an additional component of the fluid velocity visible and for a lower $\Re$ the fluid velocity is larger.  
We define the incompressibility error $e_{\DivP} = \Vert \DivP\vel \Vert_{L^{\infty}(L^2)}$ and the normal approximation error $e_N = \Vert  \vel\normal - \vnor \Vert_{L^{\infty}(L^2)}$. In figure \ref{fig:vis1} a convergence study is done which confirms third order convergence with respect to both errors. The results of $e_{\DivP}$ are similar to the results of \cite{Krause22,art:BKV23,reuther2020numerical,Brandneretal_SIAMJSC_2022} and according to \cite{hardering2021tangential} third order convergence of $e_N$ is optimal by using the surface normal $\normal$ in the penalization term. 

\subsection{Fluid like virtual clothing}
The physical coupling between normal deformation and tangential flow provides an attractive approach for interactive patterns in virtual fashion. Any body movement leads to deformation of the clothing which translates into "tangential flow", which can be visualized as interactive patterns. A realization of this idea requires dressing virtual avatars and animating them with high quality. Highly realistic physical simulation of clothing on human bodies in motion is a challenging task. We therefore do not address this approach but instead use a data-driven approach capturing dynamic clothing on humans from 4D scans, see CAPE project \cite{Pons-Moll:Siggraph2017}. CAPE is a Graph-CNN based generative model for dressing 3D meshes of human body. It is compatible with the popular body model, SMPL, and can generalize to diverse body shapes and body poses. The CAPE Dataset \cite{CAPE:CVPR:20} provides SMPL mesh registration of 4D scans of people in clothing, along with registered scans of the ground truth body shapes under clothing. 

To use this data set requires various preprocessing steps. The provided data set $\param^n_0$ does not fulfill the global area constraint within the required accuracy. We therefor modify the data set using an approach similar to eq. \eqref{eq:meshmovementII}. The second necessity for preprocessing results from the linear nature of the data set. SMPL and also almost all other popular body models use flat triangles, which is not sufficient for the numerical solution of the surface Navier-Stokes equations. The considered Taylor-Hood element with $k=3$ requires a surface mesh with at least order $k = 3$ to guarantee convergence \cite{hardering2021tangential}. The construction of higher order surfaces meshes is realized by local Gaussian smoothing. Details on these preprocessing steps are explained in \cite{kunze2022}. 

With the modified data we define $\vnor^n=(\param_k^n - \param_k^{n-1})/ \tau$ and compute the induced fluid flow $\vel$ on the hole surface $\surf_h$ by solving eq. \eqref{eq:navierstokesIV}. 
Figure \ref{fig:clothing} shows an example of such a dressed 3D mesh in motion with the induced tangential flow. As in the previous examples the simulation is done with the polynomial order of the surface and the Taylor-Hood element $k=3$. To only visualize the flow field on the clothing we use a classifier $\phi$ to define the subdomain $\surf_{cl} = \{ \coord\in\surf \vert \phi(\coord)>\delta \}$ where $\delta$ is a specific threshold. The construction of such classifiers is explained in \cite{CAPE:CVPR:20}. In figure \ref{fig:clothing} the flow field is simply visualized only on $\surf_{cl}$.

While this example only demonstrates the principle applicability of virtual interactive clothing, a realization of this concept in the booming market of virtual fashion requires further developments. 
 
\section{Conclusions}
\label{sec:conclusions}

Based on numerical approaches for fluid deformable surfaces we propose a surface finite elements method to solve the surface Navier-Stokes equations on evolving surfaces with a prescribed evolution in normal direction. Similar to the penalization of the normal component in approaches to solve vector- and tensor-valued surface partial differential equations \cite{Nestleretal_JNS_2018,jankuhn2018,Nitschkeetal_PRSA_2018,Reutheretal_PF_2018,nestler2019finite,Hansboetal_IMAJNA_2020,reuther2020numerical,Brandneretal_SIAMJSC_2022,Bachinietal_arXiv_2022} we penalize the deformation of the surface in normal direction. In contrast with previous approaches \cite{Reutheretal_MMS_2015,Reutheretal_MMS_2018,Nitschkeetal_PRF_2019,olshanskii2022tangential,olshanskii2023} treating the surface Navier-Stokes equations for the tangential component we believe the approach to be advantageous as the formulation is simpler and we expect analytical results for stationary surfaces to be applicable. The numerical results at least demonstrate the same convergence properties. However, a detailed comparison of the different approaches and an extension of the numerical analysis remains to be done.   

Interesting findings are the constraint on the prescribed normal velocity to be globally area conserving, see eq. \eqref{eq:vnorconst} and the result in figure \ref{fig:vis1} that the evolution is independent on the Reynolds number $\Re$.

The robustness of the approach allows us to apply it on more advanced surface meshes. We demonstrate the principle applicability of the approach on high-resolution 3D scans of clothed bodies in motion from the CAPE project \cite{Pons-Moll:Siggraph2017,CAPE:CVPR:20} to generate interactive virtual fluid-like clothes. The approach builds on a multi-cloth 3D model of the body and clothing. While the numerical solution requires various preprocessing steps, the approach is general and also can be applied to other body models or physical simulations of clothing. In principle, this paves the way for virtual fashion design going beyond rebuilding realistic garments in the virtual world. 

\begin{acknowledgement}
  We acknowledge computing resources provided by ZIH at TU Dresden and by JSC at FZ J\"ulich, within projects WIR and PFAMDIS, respectively. 
  This work was supported by the German Research Foundation (DFG) within the Research Unit ``Vector- and Tensor-Valued Surface PDEs'' (FOR 3013).
  Data that support the findings of this study are available from the corresponding author upon reasonable request. We further acknowledge the CAPE project for providing their 3D data of clothed bodies in motion.
\end{acknowledgement}

\bibliographystyle{pamm}
\bibliography{lib}     

\end{document}